\documentclass[a4 paper,10 pt,reqno]{amsart}
\usepackage{amsmath,graphics,hyperref}
\usepackage{amsfonts}
\usepackage{amsthm}
\usepackage{mathrsfs}  
\newtheorem{theorem}{Theorem}[section]
\newtheorem{lemma}[theorem]{Lemma}
\newtheorem{corollary}[theorem]{Corollary}
\newtheorem{remark}[theorem]{Remark}
\numberwithin{equation}{section} \makeatletter
\@namedef{subjclassname@2010}{ 2010 Mathematics Subject
Classification} \makeatother
\title{ Integral mean estimates for the polar derivative of a polynomial}
\author{N. A. Rather and Suhail Gulzar}
 \address{Department of Mathematics \\
    University of Kashmir \\
   Srinagar, Hazratbal 190006
   \\ India}
 \email{dr.narather@gmail.com}
 \email{sgmattoo@gmail.com}
\date{}

\begin{document}
\maketitle
\begin{abstract}
 \indent Let $ P(z) $ be a polynomial of degree $ n $ having all zeros in $|z|\leq k$ where $k\leq 1,$ then it was proved by Dewan \textit{et al} \cite{d} that for every real or complex number $\alpha$ with $|\alpha|\geq k$  and each $r\geq 0$
 $$   n(|\alpha|-k)\left\{\int\limits_{0}^{2\pi}\left|P\left(e^{i\theta}\right)\right|^r d\theta\right\}^{\frac{1}{r}}\leq\left\{ \int\limits_{0}^{2\pi}\left|1+ke^{i\theta}\right|^r d\theta\right\}^{\frac{1}{r}}\underset{|z|=1}{Max}|D_\alpha P(z)|. $$
  \indent In this paper, we shall present a refinement and generalization of above result and also extend it to the class of polynomials $P(z)=a_nz^n+\sum_{\nu=\mu}^{n}a_{n-\nu}z^{n-\nu},$ $1\leq\mu\leq n,$ having all its zeros in $|z|\leq k$ where $k\leq 1$ and thereby obtain certain generalizations of above and many other known results.
\end{abstract}

\section{\bf Introduction and statement of results}

Let $P(z)$ be a polynomial of degree $n$. It was shown by Tur\'{a}n \cite{t} that if $P(z)$ has all its zeros in $|z|\leq 1,$ then
\begin{equation}\label{1}
n\underset{\left|z\right|=1}{Max}\left|P(z)\right|\leq  2\underset{\left|z\right|=1}{Max}\left|P^{\prime}(z)\right|.
\end{equation}
Inequality \eqref{1} is best possible with equality holds for $P(z)=\alpha z^n+\beta$ where $|\alpha|=|\beta|.$ The above inequality \eqref{1} of Tur\'{a}n \cite{t} was generalized by Malik \cite{m69}, who proved that if $P(z)$ is a polynomial of degree $n$ having all its zeros in $|z| \leq k,$ where $k\leq 1$, then
\begin{equation}\label{2}
\underset{\left|z\right|=1}{Max}\left|P^{\prime}(z)\right|\geq        \frac{n}{1+k}\underset{\left|z\right|=1}{Max}\left|P(z)\right|.
\end{equation} 
where as for $k\geq 1,$ Govil \cite{g73} showed that
  \begin{equation}\label{3}
\underset{\left|z\right|=1}{Max}\left|P^{\prime}(z)\right|\geq        \frac{n}{1+k^n}\underset{\left|z\right|=1}{Max}\left|P(z)\right|,
\end{equation}
Both the above inequalities \eqref{2} and \eqref{3} are best possible, with equality in \eqref{2} holding for $P(z)=(z +k)^n$, where $k\geq 1.$ While in \eqref{3} the equality holds for the polynomial  $P(z)=\alpha z^n+\beta k^n$ where $|\alpha|=|\beta|.$ \\
\indent As a refinement of \eqref{2}, Aziz and Shah \cite{aw} proved if $P(z)$ is a polynomial of degree $n$ having all its zeros in $|z| \leq k,$ where $k\leq 1$, then
\begin{align}\label{3'}
\underset{\left|z\right|=1}{Max}\left|P^{\prime}(z)\right|\geq        \frac{n}{1+k}\left\{\underset{\left|z\right|=1}{Max}\left|P(z)\right|+\dfrac{1}{k^{n-1}}\underset{|z|=1}{Min}|P(z)|\right\}.
\end{align} 
\indent Let $D_\alpha P(z)$ denotes the polar derivative of the polynomial $P(z)$ of degree $n$ with respect to the point $\alpha,$ then
$$D_\alpha P(z)=nP(z)+(\alpha-z)P^{\prime}(z). $$
The polynomial $D_\alpha P(z)$ is a polynomial of degree at most $n-1$ and it generalizes the ordinary derivative in the sense that
$$\underset{\alpha\rightarrow\infty}{Lim}\left[\dfrac{D_\alpha P(z)}{\alpha}\right]=P^{\prime}(z). $$
\indent Aziz and Rather \cite{ar98} extends \eqref{2} to polar derivatives of a polynomial and proved that if all the zeros of $P(z)$ lie in $|z|\leq k$ where $k\leq 1$ then for every real or complex number $\alpha$ with $|\alpha|\geq k,$
\begin{equation}\label{4}
\underset{\left|z\right|=1}{Max}\left|D_\alpha P(z)\right|\geq n\left(\dfrac{|\alpha|-k}{1+k}\right)\underset{\left|z\right|=1}{Max}\left|P(z)\right|.
\end{equation}
 \indent For the class of polynomials $P(z)=a_nz^n+\sum_{\nu=\mu}^{n}a_{n-\nu}z^{n-\nu},$ $1\leq\mu\leq n,$  of degree $n$ having all its zeros in $|z|\leq k$ where $k\leq 1,$ Aziz and Rather \cite{ar03} proved that if $\alpha $ is real or complex number with $|\alpha|\geq k^\mu$ then
 \begin{equation}\label{4'}
\underset{\left|z\right|=1}{Max}\left|D_\alpha P(z)\right|\geq n\left(\dfrac{|\alpha|-k^\mu}{1+k^\mu}\right)\underset{\left|z\right|=1}{Max}\left|P(z)\right|.
 \end{equation}
\indent Malik \cite{m84} obtained a generalization of \eqref{1} in the sense that the left-hand side
of \eqref{1} is replaced by a factor involving the integral mean of $|P(z)|$ on $|z|=1.$ In fact he proved that if $P(z)$ has all its zeros in $|z|\leq 1,$ then for each $q>0,$
 \begin{equation}\label{5}
 n\left\{\int\limits_{0}^{2\pi}\left|P\left(e^{i\theta}\right)\right|^q d\theta\right\}^{1/q}\leq\left\{ \int\limits_{0}^{2\pi}\left|1+e^{i\theta}\right|^q d\theta\right\}^{1/q}\underset{|z|=1}{Max}|P^{\prime}(z)|.
 \end{equation}
 If we let $q$ tend to infinity in \eqref{5}, we get \eqref{1}.\\
\indent The corresponding generalization of \eqref{2} which is an extension of \eqref{5} was
 obtained by Aziz \cite{a83} by proving that if $P(z)$ is a polynomial of degree $n$ having all its zeros in $|z|\leq $ where $k\geq 1,$ then for each $q\geq 1$ 
  \begin{equation}\label{6}
  n\left\{\int\limits_{0}^{2\pi}\left|P\left(e^{i\theta}\right)\right|^q d\theta\right\}^{1/q}\leq\left\{ \int\limits_{0}^{2\pi}\left|1+k^ne^{i\theta}\right|^q d\theta\right\}^{1/q}\underset{|z|=1}{Max}|P^{\prime}(z)|.
  \end{equation}
  The result is best possible and equality in \eqref{4} holds for the polynomial $P(z)=\alpha z^n+\beta k^n$ where $|\alpha|=|\beta|.$\\
 \indent As a generalization of inequality \ref{4}, Dewan \textit{et al} \cite{d} obtained an $L^p$ inequality for the polar derivative of a polynomial and proved the following:\\
 
 \begin{theorem}\label{ta}
 If $P(z)$ is a polynomial of degree $ n $ having all its zeros in $|z|\leq k,$ where $k\leq 1,$ then for every real or complex number $\alpha$ with $|\alpha|\geq k$ and for each $r>0,$
 \begin{equation}\label{7}
  n(|\alpha|-k)\left\{\int\limits_{0}^{2\pi}\left|P\left(e^{i\theta}\right)\right|^r d\theta\right\}^{\frac{1}{r}}\leq\left\{ \int\limits_{0}^{2\pi}\left|1+ke^{i\theta}\right|^r d\theta\right\}^{\frac{1}{r}}\underset{|z|=1}{Max}|D_\alpha P(z)|.
 \end{equation} 
 \end{theorem}
 
  In this paper, we consider the class of polynomials $P(z)=a_nz^n+\sum_{j=\mu}^{n}a_{n-j}z^{n-j},$ $1\leq\mu\leq n,$ having all its zeros in $|z|\leq k$ where $k\leq 1$ and establish some improvements and generalizations of inequalities \eqref{1},\eqref{2},\eqref{4},\eqref{6} and \eqref{7}.\\
  
  \indent In this direction, we first present the following interesting results which yields \eqref{7} as a special case.
  \begin{theorem}\label{t1'}
If $P(z)$ is a polynomial of degree $ n $ having all its zeros in $|z|\leq k,$ where $k\leq 1,$ then for every real or complex $\alpha ,$ $\beta$ with $|\alpha|\geq k,$ $|\beta|\leq 1$ and for each $r>0,$ $p>1,$ $q>1$ with $p^{-1}+q^{-1}=1,$ we have
 \begin{equation}\label{te1'}
n(|\alpha|-k)\left\{\int\limits_{0}^{2\pi}\left|P(e^{i\theta})+\beta\dfrac{m}{k^{n-1}}\right|^rd\theta\right\}^{\frac{1}{r}}\leq\left\{\int\limits_{0}^{2\pi}|1+k e^{i\theta}|^{pr}d\theta\right\}^{\frac{1}{pr}}\left\{\int\limits_{0}^{2\pi}|D_\alpha P(e^{i\theta})|^{qr}d\theta\right\}^{\frac{1}{qr}}
 \end{equation}
 where $m={Min}_{|z|=k}|P(z)|.$
  \end{theorem}
 If we take $\beta=0,$ we get the following result.
 \begin{corollary}\label{c2'}
If $P(z)$ is a polynomial of degree $ n $ having all its zeros in $|z|\leq k,$ where $k\leq 1,$ then for every real or complex $\alpha ,$  with $|\alpha|\geq k$ and for each $r>0,$ $p>1,$ $q>1$ with $p^{-1}+q^{-1}=1,$ we have
\begin{equation}\label{ce2'}
n(|\alpha|-k)\left\{\int\limits_{0}^{2\pi}\left|P(e^{i\theta})\right|^rd\theta\right\}^{\frac{1}{r}}\leq\left\{\int\limits_{0}^{2\pi}|1+k e^{i\theta}|^{pr}d\theta\right\}^{\frac{1}{pr}}\left\{\int\limits_{0}^{2\pi}|D_\alpha P(e^{i\theta})|^{qr}d\theta\right\}^{\frac{1}{qr}}.
\end{equation}
 \end{corollary} 
  \begin{remark}
\textnormal{Theorem \ref{ta} follows from \eqref{ce2'} by letting  $q\rightarrow\infty$ (so that $p\rightarrow 1$) in Corollary \ref{c2'}. If we divide both sides of inequality \eqref{ce2'} by $|\alpha|$ and make $\alpha\rightarrow\infty,$ we get \eqref{4}}. 
  \end{remark}
  Dividing the two sides of \eqref{te1'} by $|\alpha|$ and letting $|\alpha|\rightarrow\infty$, we get the following result.
  \begin{corollary}
  If $P(z)$ is a polynomial of degree $ n $ having all its zeros in $|z|\leq k,$ where $k\leq 1,$ then for every real or complex  $\beta$ with $|\beta|\leq 1$ and for each $r>0,$ $p>1,$ $q>1$ with $p^{-1}+q^{-1}=1,$ we have
   \begin{equation}\label{z}
  n\left\{\int\limits_{0}^{2\pi}\left|P(e^{i\theta})+\beta\dfrac{m}{k^{n-1}}\right|^rd\theta\right\}^{\frac{1}{r}}\leq\left\{\int\limits_{0}^{2\pi}|1+k e^{i\theta}|^{pr}d\theta\right\}^{\frac{1}{pr}}\left\{\int\limits_{0}^{2\pi}|P^{\prime}(e^{i\theta})|^{qr}d\theta\right\}^{\frac{1}{qr}}
   \end{equation}
 where $m={Min}_{|z|=k}|P(z)|.$
  \end{corollary}
 If we let $q\rightarrow\infty$ in \eqref{z}, we get the following corollary. 
  \begin{corollary}
   If $P(z)$ is a polynomial of degree $ n $ having all its zeros in $|z|\leq k,$ where $k\leq 1,$ then for every real or complex  $\beta$ with $|\beta|\leq 1$ and for each $r>0,$ we have
    \begin{equation}\label{x}
   n\left\{\int\limits_{0}^{2\pi}\left|P(e^{i\theta})+\beta\dfrac{m}{k^{n-1}}\right|^rd\theta\right\}^{\frac{1}{r}}\leq\left\{\int\limits_{0}^{2\pi}|1+k e^{i\theta}|^{r}d\theta\right\}^{\frac{1}{r}}\underset{|z|=1}{Max}|P^{\prime}(z)|,
    \end{equation}
 where $m={Min}_{|z|=k}|P(z)|.$
   \end{corollary}
   \begin{remark}
   \textnormal{  If we let $r\rightarrow\infty$ in \eqref{x} and choosing argument of $\beta$ suitably with $|\beta|=1,$ we obtain \eqref{3'}. }
   \end{remark}
Next, we extend \eqref{7} to the class of polynomials $P(z)=a_nz^n+\sum_{\nu=\mu}^{n}a_{n-\nu}z^{n-\nu},$ $1\leq\mu\leq n,$ having all its zeros in $|z|\leq k,$ $k\leq 1$ and thereby obtain the following result.
  \begin{theorem}\label{t1}
If $P(z)=a_nz^n+\sum_{\nu=\mu}^{n}a_{n-\nu}z^{n-\nu},$ $1\leq\mu\leq n,$  is a polynomial of degree $n$ having all its zeros in $|z|\leq k$ where $k\leq 1,$ then for every real or complex $\alpha $ with $|\alpha|\geq k^\mu$ and for each $r>0,$ $p>1,$ $q>1$ with $p^{-1}+q^{-1}=1,$ we have
\begin{equation}\label{te1}
n(|\alpha|-k^\mu)\left\{\int\limits_{0}^{2\pi}|P(e^{i\theta})|^rd\theta\right\}^{\frac{1}{r}}\leq\left\{\int\limits_{0}^{2\pi}|1+k^\mu e^{i\theta}|^{pr}d\theta\right\}^{\frac{1}{pr}}\left\{\int\limits_{0}^{2\pi}|D_\alpha P(e^{i\theta})|^{qr}d\theta\right\}^{\frac{1}{qr}}.
\end{equation}
  \end{theorem}
  \begin{remark}
  \textnormal{ We let $r\rightarrow \infty$ and $p\rightarrow\infty$ (so that $q\rightarrow 1$) in \eqref{te1}, we get inequality \eqref{4'}.}
\end{remark}
If we divide both sides of \eqref{te1} by $|\alpha|$ and make $\alpha\rightarrow\infty,$ we get the following result.
\begin{corollary}\label{c1}
If $P(z)=a_nz^n+\sum_{\nu=\mu}^{n}a_{n-\nu}z^{n-\nu},$ $1\leq\mu\leq n,$  is a polynomial of degree $n$ having all its zeros in $|z|\leq k$ where $k\leq 1,$ then for  for each $r>0,$ $p>1,$ $q>1$ with $p^{-1}+q^{-1}=1,$ we have
\begin{equation}\label{ce1}
n\left\{\int\limits_{0}^{2\pi}|P(e^{i\theta})|^rd\theta\right\}^{\frac{1}{r}}\leq\left\{\int\limits_{0}^{2\pi}|1+k^\mu e^{i\theta}|^{pr}d\theta\right\}^{\frac{1}{pr}}\left\{\int\limits_{0}^{2\pi}|P^{\prime}(e^{i\theta})|^{qr}d\theta\right\}^{\frac{1}{qr}}.
\end{equation}
\end{corollary}
 Letting $q\rightarrow\infty$ (so that $p\rightarrow 1$) in \eqref{te1}, we get the following result:
\begin{corollary}\label{c3}
If $P(z)=a_nz^n+\sum_{\nu=\mu}^{n}a_{n-\nu}z^{n-\nu},$ $1\leq\mu\leq n,$ where $1\leq \mu\leq n,$ is a polynomial of degree $ n $ having all its zeros in $|z|\leq k,$ where $k\leq 1,$ then for every real or complex number $\alpha$ with $|\alpha|\geq k^\mu$ and for each $r>0,$
 \begin{equation}\label{ce3}
  n(|\alpha|-k^\mu)\left\{\int\limits_{0}^{2\pi}\left|P\left(e^{i\theta}\right)\right|^r d\theta\right\}^{\frac{1}{r}}\leq\left\{ \int\limits_{0}^{2\pi}\left|1+k^\mu e^{i\theta}\right|^r d\theta\right\}^{\frac{1}{r}}\underset{|z|=1}{Max}|D_\alpha P(z)|.
 \end{equation} 
\end{corollary} 

As a generalization of Theorem \ref{t1}, we present the following result:
\begin{theorem}\label{t2}
If $P(z)=a_nz^n+\sum_{\nu=\mu}^{n}a_{n-\nu}z^{n-\nu}$ where $1\leq \mu\leq n,$ is a polynomial of degree $n$ having all its zeros in $|z|\leq k$ where $k\leq 1,$ then for every real or complex $\alpha $ with $|\alpha|\geq k^\mu$ and for each $r>0,$ $p>1,$ $q>1$ with $p^{-1}+q^{-1}=1,$ we have
\begin{align}\label{te2}
n(|\alpha|-k^\mu)\left\{\int\limits_{0}^{2\pi}|P(e^{i\theta})+\beta m|^rd\theta\right\}^{\frac{1}{r}}\leq\left\{\int\limits_{0}^{2\pi}|1+k^\mu e^{i\theta}|^{pr}d\theta\right\}^{\frac{1}{pr}}\left\{\int\limits_{0}^{2\pi}|D_\alpha P(e^{i\theta})|^{qr}d\theta\right\}^{\frac{1}{qr}}
\end{align}
 where $m={Min}_{|z|=k}|P(z)|.$
\end{theorem}

If we divide both sides by $|\alpha|$ and make $\alpha\rightarrow\infty,$ we get the following result:
\begin{corollary}\label{c4}
If $P(z)=a_nz^n+\sum_{\nu=\mu}^{n}a_{n-\nu}z^{n-\nu},$ $1\leq\mu\leq n,$  is a polynomial of degree $n$ having all its zeros in $|z|\leq k$ where $k\leq 1,$ then for  for each $r>0,$ $p>1,$ $q>1$ with $p^{-1}+q^{-1}=1,$ we have
\begin{equation}\label{ce4}
n\left\{\int\limits_{0}^{2\pi}|P(e^{i\theta})+\beta m|^rd\theta\right\}^{\frac{1}{r}}\leq\left\{\int\limits_{0}^{2\pi}|1+k^\mu e^{i\theta}|^{pr}d\theta\right\}^{\frac{1}{pr}}\left\{\int\limits_{0}^{2\pi}|P^{\prime}(e^{i\theta})|^{qr}d\theta\right\}^{\frac{1}{qr}}
\end{equation}
 where $m={Min}_{|z|=k}|P(z)|.$
\end{corollary}
 Letting $q\rightarrow\infty$ (so that $p\rightarrow 1$) in \eqref{te1}, we get the following result:
\begin{corollary}\label{c6}
If $P(z)=a_nz^n+\sum_{\nu=\mu}^{n}a_{n-\nu}z^{n-\nu}$  where $1\leq \mu\leq n,$ is a polynomial of degree $ n $ having all its zeros in $|z|\leq k,$ where $k\leq 1,$ then for every real or complex number $\alpha$ with $|\alpha|\geq k^\mu$ and for each $r>0,$
 \begin{equation}\label{ce6}
  n(|\alpha|-k^\mu)\left\{\int\limits_{0}^{2\pi}\left|P\left(e^{i\theta}\right)+\beta m\right|^r d\theta\right\}^{\frac{1}{r}}\leq\left\{ \int\limits_{0}^{2\pi}\left|1+k^\mu e^{i\theta}\right|^r d\theta\right\}^{\frac{1}{r}}\underset{|z|=1}{Max}|D_\alpha P(z)|
 \end{equation} 
 where $m={Min}_{|z|=k}|P(z)|.$
\end{corollary}

\section{\bf Lemmas}

For the proofs of the theorems, we need the following Lemmas:
\begin{lemma}\label{l1'}
If $P(z)$ is a polynomial of degree almost $n$ having all its zeros in in $|z|\leq k$ $k\leq 1$ then for $|z|=1,$
\begin{equation}
|Q^{\prime}(z)|+\dfrac{nm}{k^{n-1}}\leq k|P^{\prime}(z)|,
\end{equation}
where $Q(z)=z^n\overline{P(1/\overline{z})}$ and $m={Min}_{|z|=k}|P(z)|.$
\end{lemma}
The above Lemma is due to Govil and McTume \cite{gm}.
\begin{lemma}\label{l1}
Let $P(z)=a_0+\sum_{\nu=\mu}^{n}a_{\nu}z^{\nu},$ $1\leq\mu\leq n,$ is a polynomial of degree $n,$ which does not vanish for $|z|<k,$ where $k\geq 1$  then for $|z|=1,$ 
\begin{equation}\label{le1}
k^\mu|P^{\prime}(z)|\leq |Q^{\prime}(z)|,
\end{equation}
where $Q(z)=z^n\overline{P(1/\overline{z})}.$
\end{lemma}
The above Lemma is due to Chan and Malik \cite{cm}. By applying Lemma \ref{l1} to the polynomial $ z^n\overline{P(1/\overline{z})} ,$ one can easily deduce:
\begin{lemma}\label{l2}
Let $P(z)=a_nz^n+\sum_{\nu=\mu}^{n}a_{n-\nu}z^{n-\nu},$ $1\leq\mu\leq n,$ is a polynomial of degree $n,$ having all its zeros in $|z|\leq k,$ where $k\leq 1$  then for $|z|=1$ 
\begin{equation}
k^\mu|P^{\prime}(z)|\geq |Q^{\prime}(z)|,
\end{equation}
where $Q(z)=z^n\overline{P(1/\overline{z})}.$
\end{lemma}

\section{\bf Proof of Theorems}

\begin{proof}[\bf Proof of Theorem \ref{t1'}]
 Let $Q(z)=z^n\overline{P(1/\overline{z})}$ then $P(z)=z^n\overline{Q(1/\overline{z})}$ and it can be easily verified that for $|z|=1,$
\begin{equation}\label{p1}
|Q^{\prime}(z)|= |nP(z)-zP^{\prime}(z)|\,\,\,\,  \textnormal{and}\,\,\,\,  |P^{\prime}(z)|= |nQ(z)-zQ^{\prime}(z)|.
\end{equation}
By Lemma \eqref{l1'}, we have for every $\beta $ with $|\beta|\leq 1$ and  $|z|=1,$
\begin{align}\label{p2}
\left|Q^{\prime}(z)+\bar{\beta}\dfrac{ nmz^{n-1}}{k^{n-1}}\right|\leq|Q^{\prime}(z)|+\dfrac{nm}{k^{n-1}}\leq k|P^{\prime}(z)|.
\end{align}
Using \eqref{p1} in \eqref{p2}, for $|z|=1$ we have
\begin{equation}\label{p3}
\left|Q^{\prime}(z)+\bar{\beta}\dfrac{ nmz^{n-1}}{k^{n-1}}\right|\leq k|nQ(z)-zQ^{\prime}(z)|.
\end{equation}
By Lemma \ref{l2} with $\mu=1,$ for every real or complex number $\alpha$ with $|\alpha|\geq k$ and $|z|=1,$ we have
\begin{align}\label{p6}\nonumber
|D_\alpha P(z)|&\geq |\alpha||P^{\prime}(z)|-|Q^{\prime}(z)|\\&\geq (|\alpha|-k)|P^{\prime}(z)|.
\end{align}
Since $P(z)$ has all its zeros in $|z|\leq k\leq 1,$ it follows by Gauss-Lucas Theorem that all the zeros of $P^{\prime}(z)$ also lie in $|z|\leq k\leq 1.$ This implies that the polynomial
\begin{equation*}
z^{n-1}\overline{P^{\prime}(1/\overline{z})}\equiv nQ(z)-zQ^{\prime}(z)
\end{equation*}
does not vanish in $|z|<1.$ Therefore, it follows from \eqref{p3} that the function
\begin{equation*}
w(z)=\dfrac{z\left(Q^{\prime}(z)+\bar{\beta}\dfrac{nmz^{n-1}}{k^{n-1}}\right)}{k\left(nQ(z)-zQ^{\prime}(z)\right)}
\end{equation*}
is analytic for $|z|\leq 1$ and $|w(z)|\leq 1$ for $|z|=1.$ Furthermore, $w(0)=0.$ Thus the function $1+k w(z)$ is subordinate to the function $1+k z$ for $|z|\leq 1.$ Hence by a well known property of subordination \cite{h}, we have
\begin{equation}\label{p4}
\int\limits_{0}^{2\pi}\left|1+k w(e^{i\theta})\right|^rd\theta\leq\int\limits_{0}^{2\pi}\left|1+k e^{i\theta}\right|^rd\theta,\,\,\,r>0.
\end{equation}
Now
\begin{equation*}
1+k w(z)=\dfrac{n\left(Q(z)+\bar{\beta}\dfrac{mz^{n}}{k^{n-1}}\right)}{nQ(z)-zQ^{\prime}(z)},
\end{equation*} 
and
\begin{equation*}
|P^{\prime}(z)|=|z^{n-1}\overline{P^{\prime}(1/\overline{z})}|=|nQ(z)-zQ^{\prime}(z)|,\,\,\,\textrm{for}\,\,\,|z|=1,
\end{equation*}
therefore for $|z|=1,$
\begin{equation*}
n\left|Q(z)+\bar{\beta}\dfrac{mz^{n}}{k^{n-1}}\right|=|1+ kw(z)||nQ(z)-zQ^{\prime}(z)|=|1+k w(z)||P^{\prime}(z)|.
\end{equation*}
equivalently,
\begin{align*}
n\left|z^n\overline{P(1/\overline{z})}+\bar{\beta}\dfrac{mz^{n}}{k^{n-1}}\right|=|1+k w(z)||P^{\prime}(z)|.
\end{align*}
This implies 
\begin{equation}\label{p5}
n\left|P(z)+\beta\dfrac{m}{k^{n-1}}\right|=|1+k w(z)||P^{\prime}(z)|\,\,\,\,\textnormal{for}\,\,\,\,|z|=1.
\end{equation}
From \eqref{p6} and \eqref{p5}, we deduce that for $r>0,$
\begin{equation*}
n^r(|\alpha|-k)^r\int\limits_{0}^{2\pi}\left|P(e^{i\theta})+\beta\dfrac{m}{k^{n-1}}\right|^rd\theta\leq\int\limits_{0}^{2\pi}|1+ k w(e^{i\theta})|^r|D_\alpha P(e^{i\theta})|^rd\theta.
\end{equation*}
This gives with the help of H\"{o}lder's inequality and using \eqref{p4}, for $p>1,$ $q>1$ with $p^{-1}+q^{-1}=1,$
\begin{equation*}
n^r(|\alpha|-k)^r\int\limits_{0}^{2\pi}\left|P(e^{i\theta})+\beta\dfrac{m}{k^{n-1}}\right|^rd\theta\leq\left(\int\limits_{0}^{2\pi}|1+k e^{i\theta}|^{pr}d\theta\right)^{1/p}\left(\int\limits_{0}^{2\pi}|D_\alpha P(e^{i\theta})|^{qr}d\theta\right)^{1/q},
\end{equation*}
equivalently,
\begin{equation*}
n(|\alpha|-k^\mu)\left\{\int\limits_{0}^{2\pi}\left|P(e^{i\theta})+\beta\dfrac{m}{k^{n-1}}\right|^rd\theta\right\}^{\frac{1}{r}}\leq\left\{\int\limits_{0}^{2\pi}|1+k e^{i\theta}|^{pr}d\theta\right\}^{\frac{1}{pr}}\left\{\int\limits_{0}^{2\pi}|D_\alpha P(e^{i\theta})|^{qr}d\theta\right\}^{\frac{1}{qr}}
\end{equation*}
which proves the desired result.
\end{proof}
\begin{proof}[\bf Proof of Theorem \ref{t1}]
Since $P(z)$ has all its zeros in $|z|\leq k,$ therefore, by using Lemma \ref{l2} we have for $|z|=1,$
\begin{equation}\label{t1p2}
|Q^{\prime}(z)|\leq k^\mu |nQ(z)-zQ^{\prime}(z)|.
\end{equation}
Now for every real or complex number $\alpha$ with $|\alpha|\geq k^\mu,$  we have
\begin{align*}
|D_\alpha P(z)|&=|nP(z)+(\alpha-z)P^{\prime}(z)|\\&\geq |\alpha||P^{\prime}(z)|-|nP(z)-zP^{\prime}(z)|,
\end{align*} 
by using \eqref{p1} and Lemma \ref{l2}, for $|z|=1,$  we get
\begin{align}\nonumber\label{t1p3}
|D_\alpha P(z)|&\geq |\alpha||P^{\prime}(z)|-|Q^{\prime}(z)|\\&\geq (|\alpha|-k^\mu)|P^{\prime}(z)|.
\end{align}
Since $P(z)$ has all its zeros in $|z|\leq k\leq 1,$ it follows by Gauss-Lucas Theorem that all the zeros of $P^{\prime}(z)$ also lie in $|z|\leq k\leq 1.$ This implies that the polynomial
\begin{equation*}
z^{n-1}\overline{P^{\prime}(1/\overline{z})}\equiv nQ(z)-zQ^{\prime}(z)
\end{equation*}
does not vanish in $|z|<1.$ Therefore, it follows from \eqref{t1p2} that the function
\begin{equation*}
w(z)=\dfrac{zQ^{\prime}(z)}{k^\mu\left(nQ(z)-zQ^{\prime}(z)\right)}
\end{equation*}
is analytic for $|z|\leq 1$ and $|w(z)|\leq 1$ for $|z|=1.$ Furthermore, $w(0)=0.$ Thus the function $1+k^\mu w(z)$ is subordinate to the function $1+k^\mu z$ for $|z|\leq 1.$ Hence by a well known property of subordination \cite{h}, we have
\begin{equation}\label{t1p4}
\int\limits_{0}^{2\pi}\left|1+k^\mu w(e^{i\theta})\right|^rd\theta\leq\int\limits_{0}^{2\pi}\left|1+k^\mu e^{i\theta}\right|^rd\theta,\,\,\,r>0.
\end{equation}
Now
\begin{equation*}
1+k^\mu w(z)=\dfrac{nQ(z)}{nQ(z)-zQ^{\prime}(z)},
\end{equation*} 
and
\begin{equation*}
|P^{\prime}(z)|=|z^{n-1}\overline{P^{\prime}(1/\overline{z})}|=|nQ(z)-zQ^{\prime}(z)|,\,\,\,\textrm{for}\,\,\,|z|=1,
\end{equation*}
therefore, for $|z|=1,$
\begin{equation}\label{t1p5}
n|Q(z)|=|1+k^\mu w(z)||nQ(z)-zQ^{\prime}(z)|=|1+k^\mu w(z)||P^{\prime}(z)|.
\end{equation}
From \eqref{t1p3} and \eqref{t1p5}, we deduce that for $r>0,$
\begin{equation*}
n^r(|\alpha|-k^\mu)^r\int\limits_{0}^{2\pi}|P(e^{i\theta})|^rd\theta\leq\int\limits_{0}^{2\pi}|1+k^\mu w(e^{i\theta})|^r|D_\alpha P(e^{i\theta})|^rd\theta.
\end{equation*}
This gives with the help of H\"{o}lder's inequality and \eqref{t1p4}, for $p>1,$ $q>1$ with $p^{-1}+q^{-1}=1,$
\begin{equation*}
n^r(|\alpha|-k^\mu)^r\int\limits_{0}^{2\pi}|P(e^{i\theta})|^rd\theta\leq\left(\int\limits_{0}^{2\pi}|1+k^\mu e^{i\theta}|^{pr}d\theta\right)^{1/p}\left(\int\limits_{0}^{2\pi}|D_\alpha P(e^{i\theta})|^{qr}d\theta\right)^{1/q},
\end{equation*}
equivalently,
\begin{equation*}
n(|\alpha|-k^\mu)\left\{\int\limits_{0}^{2\pi}|P(e^{i\theta})|^rd\theta\right\}^{\frac{1}{r}}\leq\left\{\int\limits_{0}^{2\pi}|1+k^\mu e^{i\theta}|^{pr}d\theta\right\}^{\frac{1}{pr}}\left\{\int\limits_{0}^{2\pi}|D_\alpha P(e^{i\theta})|^{qr}d\theta\right\}^{\frac{1}{qr}}
\end{equation*}
which proves the desired result.

\end{proof}

\begin{proof}[\bf Proof of Theorem \ref{t2}]
Let  $m=Min_{|z|=k}|P(z)|,$ so that $m\leq |P(z)|$ for $|z|=k.$ If $P(z)$ has a zero on $|z|=k$ then $m=0$ and result follows from Theorem \ref{t1}. Henceforth we suppose that all the zeros of $P(z)$ lie in $|z|< k.$  Therefore for every $\beta$ with $|\beta|<1,$ we have $|m\beta|<|P(z)|$ for $|z|=k.$  Since $P(z)$ has all its zeros in $|z|<k\leq 1,$ it follows by Rouche's theorem that all the zeros of $F(z)=P(z)+\beta m$ lie in $|z|< k\leq 1.$ 
If $G(z)=z^n\overline{F(1/\overline{z})}=Q(z)+\bar{\beta} mz^n,$ then by applying Lemma \ref{l2} to polynomial $F(z)=P(z)+\beta m,$ we have for $|z|=1,$
\begin{equation*}
|G^{\prime}(z)|\leq k^\mu |F^{\prime}(z)|.
\end{equation*}
This gives
\begin{equation}\label{t2p4}
|Q^{\prime}(z)+nm\bar{\beta} z^{n-1}|\leq k^\mu|P^{\prime}(z)|.
\end{equation}
Using \eqref{p1} in \eqref{t2p4}, for $|z|=1$ we have
\begin{equation}\label{t2p4'}
|Q^{\prime}(z)+nm\bar{\beta} z^{n-1}|\leq k^\mu|nQ(z)-zQ^{\prime}(z)|
\end{equation}
Since $P(z)$ has all its zeros in $|z|< k\leq 1,$ it follows by Gauss-Lucas Theorem that all the zeros of $P^{\prime}(z)$ also lie in $|z|< k\leq 1.$ This implies that the polynomial
\begin{equation*}
z^{n-1}\overline{P^{\prime}(1/\overline{z})}\equiv nQ(z)-zQ^{\prime}(z)
\end{equation*}
does not vanish in $|z|<1.$ Therefore, it follows from \eqref{t2p4'} that the function
\begin{equation*}
w(z)=\dfrac{z(Q^{\prime}(z)+nm\bar{\beta} z^{n-1})}{k^\mu\left(nQ(z)-zQ^{\prime}(z)\right)}
\end{equation*}
is analytic for $|z|\leq 1$ and $|w(z)|\leq 1$ for $|z|=1.$ Furthermore, $w(0)=0.$ Thus the function $1+k^\mu w(z)$ is subordinate to the function $1+k^\mu z$ for $|z|\leq 1.$ Hence by a well known property of subordination \cite{h}, we have
\begin{equation}\label{t2p5}
\int\limits_{0}^{2\pi}\left|1+k^\mu w(e^{i\theta})\right|^rd\theta\leq\int\limits_{0}^{2\pi}\left|1+k^\mu e^{i\theta}\right|^rd\theta,\,\,\,r>0.
\end{equation}
Now
\begin{equation*}
1+k^\mu w(z)=\dfrac{n(Q(z)+m\bar{\beta} z^{n})}{nQ(z)-zQ^{\prime}(z)},
\end{equation*} 
and
\begin{equation*}
|P^{\prime}(z)|=|z^{n-1}\overline{P^{\prime}(1/\overline{z})}|=|nQ(z)-zQ^{\prime}(z)|,\,\,\,\textrm{for}\,\,\,|z|=1,
\end{equation*}
therefore, for $|z|=1,$
\begin{equation*}
n|Q(z)+m\bar{\beta} z^{n}|=|1+k^\mu w(z)||nQ(z)-zQ^{\prime}(z)|=|1+k^\mu w(z)||P^{\prime}(z)|.
\end{equation*}
This implies 
\begin{equation}\label{t2p6}
n|G(z)|=|1+k^\mu w(z)||nQ(z)-zQ^{\prime}(z)|=|1+k^\mu w(z)||P^{\prime}(z)|.
\end{equation}
Since $|F(z)|=|G(z)|$ for $|z|=1,$ therefore, from \eqref{t2p6} we get
\begin{equation}\label{t2p7}
n|P(z)+\beta m|=|1+k^\mu w(z)||P^{\prime}(z)|\,\,\,\,\textnormal{for}\,\,\,\,\,|z|=1.
\end{equation}
From \eqref{t1p3} and \eqref{t2p7}, we deduce that for $r>0,$
\begin{equation*}
n^r(|\alpha|-k^\mu)^r\int\limits_{0}^{2\pi}|P(e^{i\theta})+\beta m|^rd\theta\leq\int\limits_{0}^{2\pi}|1+k^\mu w(e^{i\theta})|^r|D_\alpha P(e^{i\theta})|^rd\theta.
\end{equation*}
This gives with the help of H\"{o}lder's inequality in conjunction with \eqref{t2p5} for $p>1,$ $q>1$ with $p^{-1}+q^{-1}=1,$
\begin{equation*}
n^r(|\alpha|-k^\mu)^r\int\limits_{0}^{2\pi}|P(e^{i\theta})+\beta m|^rd\theta\leq\left(\int\limits_{0}^{2\pi}|1+k^\mu e^{i\theta}|^{pr}d\theta\right)^{1/p}\left(\int\limits_{0}^{2\pi}|D_\alpha P(e^{i\theta})|^{qr}d\theta\right)^{1/q},
\end{equation*}
equivalently,
\begin{equation*}
n(|\alpha|-k^\mu)\left\{\int\limits_{0}^{2\pi}|P(e^{i\theta})+\beta m|^rd\theta\right\}^{\frac{1}{r}}\leq\left\{\int\limits_{0}^{2\pi}|1+k^\mu e^{i\theta}|^{pr}d\theta\right\}^{\frac{1}{pr}}\left\{\int\limits_{0}^{2\pi}|D_\alpha P(e^{i\theta})|^{qr}d\theta\right\}^{\frac{1}{qr}}
\end{equation*}
which proves the desired result.
\end{proof}

\end{document}